\documentclass[12pt]{article}
\usepackage{hyperref}
\usepackage{amssymb, amsmath, amsfonts}
\usepackage{amsthm}
\usepackage[all]{xy}
\usepackage{tikz-cd}
\usepackage[margin=1in]{geometry}

\usepackage{bm}
\def\genfd{{\bm k}}
\def\id{{\rm id}}
\long\def\nodo#1{{}}

\newtheorem{thm}{Theorem}
\newtheorem{lem}{Lemma}
\newtheorem{prop}{Proposition}
\newtheorem{defn}{Definition}

\newtheorem{remark}{Remark}

\def\MR#1{} 
\usepackage{graphicx}

\title{Extending Doi--Takeuchi correspondence to bimodules and noncommutative covers}
\author{Zoran \v{S}koda}
 \date{}

\begin{document}

\maketitle
\begin{abstract}
Doi and Takeuchi, and independently Blattner, Cohen and Montgomery, have exhibited
a correspondence between cleft extensions of an algebra $U$ by a Hopf algebra $H$
and $H$-measurings on $U$ equipped with a compatible normalized 2-cocycle. Each cleft extension is isomorphic as an algebra, a left $U$-module and a right $H$-comodule to the crossed product algebra constructed from the measuring and the cocycle. Cleft $H$-extensions are a special case of $H$-Hopf-Galois extensions, describing noncommutative principal $H$-bundles with affine total and affine base space. Corings are used to describe noncommutative space covers, allowing a description of possibly nonaffine noncommutative spaces. We propose to describe the total space of a noncommutative $H$-principal bundle by an internal $E$-coring in $H$-comodules where $E$ is an $H$-comodule algebra which is itself $H$-Hopf-Galois. This article focuses on the cleft case.

As corings involve bimodules and morphisms of bimodules we extend the Doi--Takeuchi correspondence to bimodule $H$-comodules between cleft extensions on one side and bimodules of $H$-coinvariants equipped with bimodule $H$-measurings compatible with 2-cocycles on the other side. In fact, the correspondences for algebras and for bimodules are parts of a biequivalence of suitable bicategories defined in this article. As a corollary one obtains a correspondence between cleft $H$-extensions of corings and corings equipped with a nonabelian \v{C}ech 1-cocycle covalued in $H$.
\end{abstract}

\section{Introduction}

\subsection{Motivation and overview}

In noncommutative algebraic geometry, Hopf algebras are often viewed as a generalization of (the categorically formal duals to) affine algebraic groups. Faithfully flat Hopf-Galois extensions~\cite{BrzMaj:quantumgauge,Schn:prinHomSp} of algebras by Hopf algebras are analogues of torsors (in appropriate flat topology) for affine algebraic groups whose total and base spaces are affine as well. In noncommutative case we can glue pieces using flat localizations~\cite{skoda:ban,ros:NcSch}, ideals~\cite{brzwri} or faithfully flat maps~\cite{KontsRos}. If the pieces are affine, corings supply a formalism~\cite{brzwisCorings,brzwri}; if the pieces themselves are not affine we may need more general comonads~\cite{BohmSkoda,ros:NcSch,skoda:gmj}.

We may also be interested that the affine pieces are in some sense trivial principal bundles~\cite{BrzMaj:quantumgauge}. In commutative case, the algebra of functions on trivial principal bundle (when both the base and total space are affine) is the tensor product of the algebra of functions on the base and of the fiber. In noncommutative case the tensor product algebra seems to be too restrictive in major examples and one goes beyond to include either smash product algebras~\cite{Doi,montg} or more generally cocycled crossed product algebras~\cite{montg,doitakeuchi} with canonical structure of a comodule over the Hopf algebra. In both cases one should allow also objects which are isomorphic to these products in appropriate sense. For crossed product algebras these appear to be cleft Hopf-Galois extensions. Each class of cleft Hopf-Galois extensions is represented up to isomorphism by a cocycled crossed product algebra. This is the geometric meaning of the Doi--Takeuchi correspondence~\cite{BCM,doitakeuchi} reviewed in~\ref{ssec:crossedcleft}. The present article globalizes this correspondence to covers presented by corings generalizing an earlier sketch for the smash product case~\cite{skoda:hgcoring2011}. Technical elements of the present construction involve internalization of corings~\cite{Bohm:intcor} to the category of comodules over a Hopf algebra, a new (bi)cocycled crossed product construction for bimodules and a bicategorical setup unifying the algebra, bimodule and morphism of bimodules aspects of the correspondence.

This description of locally flat bundles is just a small part of the emerging picture fitting in the general philosophy of~\cite{skoda:gmj}. More general case of locally Hopf-Galois extensions is treated in~\cite{skoda:hgcoring} along with the proof(s) of the Schneider's descent theorem~\cite{Schn:prinHomSp} in this setup; this can be further extended for coalgebra extensions (in the sense of~\cite{BrzMaj:coalgbun,brzwisCorings} and entwining structures. The base change to refinements of covers using appropriate morphisms of corings relates the corresponding \v{C}ech cocycles; this has many cumbersome details and some choices, we delegate a full treatment to future publications. The associated bundles for Hopf-Galois extension of corings not only at the level of global section but noncommutative spaces presented by corings as well is in~\cite{skoda:ass}.

 An early part of this work was presented in Istanbul 2011 and widely discussed with T. Brzezi\'nski in Swansea 2014 where my stay was financed by G. Garkusha's grant. I thank l'IH\'ES for hosting my visit in 2019 where the bicategorical setup has been conceived. I thank the organizers of Hopf25 in Bruxelles where the mature version of this work has been presented. The {\tt tikzcd} code of the coherence diagrams~(\ref{eq:cohforCross}),(\ref{eq:PhiCoh}),(\ref{eq:PsiCoh}) in this paper was drawn by a free GenAI tool.

\subsection{Basics on crossed products and cleft extensions}
\label{ssec:crossedcleft}

All algebras, coalgebras and unadorned tensor products of modules are over a fixed ground field $\genfd$.
Given a bialgebra $H$ with comultiplication $\Delta$ and counit $\epsilon$, a right $H$-comodule algebra $E$ with coaction $\rho\colon E\to E\otimes H$ is said to be cleft if there is a morphism of right $H$-comodules $\gamma\colon H\to E$ which admits a convolution inverse $\gamma^{-1}\colon H\to E$, that is $\gamma^{-1}(h_{(1)})\gamma(h_{(2)}) = \epsilon(h)1_E = \gamma(h_{(1)})\gamma^{-1}(h_{(2)})$ for all $h\in H$. If $H$ admits an antipode it follows that $\rho(\gamma^{-1}(h))=\gamma^{-1}(h_{(2)})\otimes S h_{(1)}$. Here and elsewhere we use Sweedler notation without an explicit summation sign. Since $1_H$ is grouplike, $\gamma^{-1}(1_H)$ is the usual inverse of $\gamma(1_H)$ in $E$, hence the cleaving $\gamma$ can be multiplied
by $\gamma^{-1}(1_H)$ yielding a new cleaving $\gamma'$ satisfying unitality $\gamma'(1_H) = 1_E$. For the subalgebra of $H$-coinvariants $E^{\mathrm{co}H} = \{e\in E \mid \rho(e)=e\otimes 1\}$ one defines a $\genfd$-linear map
\begin{equation}
\triangleright = \triangleright_\gamma\colon H\otimes E^{\mathrm{co}H}\to E^{\mathrm{co}H},
\quad h\otimes u\mapsto \gamma(h_{(1)}) u \gamma^{-1}(h_{(2)}),
\end{equation}
which is a left $H$-measuring~\cite{sweedler,montg} in the sense that $h\triangleright 1 = \epsilon(h)1$ and $h\triangleright(u\cdot v) = (h_{(1)}\triangleright u)\cdot(h_{(2)}\triangleright v)$ for all $h\in H$, $u,v\in E^{\mathrm{co}H}$. If $\gamma$ is chosen unital, as we assume from now on, then also $1_H\triangleright_\gamma u = u$. In that case, if the $\genfd$-linear map
\begin{equation}\label{eq:sigmagamma}
\sigma = \sigma_\gamma\colon H\otimes H\to E^{\mathrm{co}H},
\quad h\otimes k\mapsto \gamma(h_{(1)})\gamma(k_{(1)})\gamma^{-1}(h_{(2)}k_{(2)})
\end{equation}
is trivial in the sense that $\sigma(h\otimes k) = \epsilon(h)\epsilon(k)1_E$ then $\triangleright$ is an action. In any case, $\sigma\in\operatorname{Hom}_\genfd(H\otimes H,E^{\mathrm{co}H})$ is a convolution invertible 2-cocycle,
\begin{equation}\label{eq:2coc}
    \sum [h_{(1)}\triangleright\sigma(k_{(1)},m_{(1)})]\sigma(h_{(2)},k_{(2)}m_{(2)})
    =\sum \sigma(h_{(1)},k_{(1)})\sigma(h_{(2)}k_{(2)},m),
\end{equation}
for all $h,k,m\in H$ and with notation $\sigma(h,k):=\sigma(h\otimes k)$.
Moreover, $\sigma$ is normalized, $\sigma(1,h) = \sigma(h,1) = \epsilon(h)1$, and it is compatible with $\triangleright$ in the sense that, for all $h,k\in H$ and $u\in E^{\mathrm{co}H}$,
\begin{equation}\label{eq:compalg}
h\triangleright (k\triangleright u) = \sum \sigma(h_{(1)},k_{(1)})
((h_{(2)}k_{(2)})\triangleright u) \sigma^{-1}(h_{(3)},k_{(3)}).
\end{equation}

\begin{defn}\label{def:crossedsystem}
A \emph{crossed system} (or \emph{complete cocycle datum}) for a bialgebra $H$ acting on
a unital associative $\genfd$-algebra $U$ is a triple $(U,\triangleright,\sigma)$ where
$\triangleright\colon H\otimes U\to U$ is a unital $H$-measuring and
$\sigma\colon H\otimes H\to U$ is a convolution invertible normalized 2-cocycle satisfying the compatibility condition~(\ref{eq:compalg}).
\end{defn}

If $\triangleright$ is any $H$-measuring on a unital associative algebra $U$ then
the formula
\begin{equation}
(u\otimes h)(v\otimes g) = u (h_{(1)}\triangleright v)\sigma(h_{(2)},g_{(1)})\otimes h_{(3)} g_{(2)}
\end{equation}
defines an associative multiplication on $U\otimes H$ with unit $1_U\otimes 1_H$ iff $\sigma$ is a normalized 2-cocycle compatible with $\triangleright$. With this multiplication, $U\otimes H$ is called a \emph{(cocycled) crossed product algebra}, $U\sharp_{\triangleright,\sigma} H$, sometimes denoted just $U\sharp_\sigma H$. For elements we write $u\sharp h$ for $u\otimes h\in U\sharp_\sigma H$. The underlying $\genfd$-module of $U\sharp_\sigma H$ becomes a left $U$-module by $u(v\sharp h) := (u v)\sharp h$ and a right $H$-comodule by coaction $\id_U\otimes\Delta$. If $H$ is a Hopf algebra with antipode $S\colon H\to H$ then the $H$-comodule map $H\to U\sharp H$, $h\mapsto 1\sharp h$, has a convolution inverse $h\mapsto\sigma^{-1}(Sh_{(2)},h_{(3)})\otimes Sh_{(1)}$,
hence $U\sharp_\sigma H$ is a cleft right $H$-comodule algebra with cleaving $h\mapsto 1\sharp h$.
If $U = E^{\mathrm{co}H}$ where $(E,\rho\colon e\mapsto e_{(0)}\otimes e_{(1)})$ is a cleft $H$-comodule algebra with cleaving $\gamma$, then
$E\ni e\mapsto e_{(0)}\gamma^{-1}(e_{(1)})\sharp e_{(2)}\in E^{\mathrm{co}H}\sharp_{\triangleright_\gamma,\sigma_\gamma}H$
is an isomorphism of $\genfd$-algebras, left $E^{\mathrm{co}H}$-modules and right $H$-comodules,
with inverse $u\sharp h\mapsto u\gamma(h)$. Conversely, if we start with an algebra $U$ with a normalized 2-cocycle $\sigma$ and a compatible unital $H$-measuring $\triangleright$, then the subalgebra of coinvariants $(U\sharp_\sigma H)^{\mathrm{co}H}$ equals $U\otimes \genfd 1_H$, clearly isomorphic to $U$, and the 2-cocycle $\sigma'$ and $H$-measuring $\triangleright'$ induced by its cleaving $h\mapsto 1\sharp h$ may be identified with $\sigma$ and $\triangleright$ via this isomorphism. This is the Doi--Takeuchi correspondence~\cite{BCM,doitakeuchi,Doi,montg} between crossed systems $(U,\triangleright,\sigma)$ and cleft $H$-comodule algebras $(E,\gamma)$: any cleft $H$-comodule algebra can be deconstructed into the algebra of $H$-coinvariants with a measuring and a cocycle, and then reconstructed, up to an isomorphism, as the cocycled crossed product algebra.

In this article, the bimodule actions are always denoted by lower dots $.$, while the product in the algebras will be either $\cdot$ or concatenation. For one-sided modules the action we choose either lower dot or concatenation, depending on context.

\section{Correspondence for bimodules}

\subsection{Bimodule measurings versus bimodules with $H$-coaction}
\label{ssec:DTbim}
\begin{defn}
  Given two $\genfd$-algebras $U$ and $V$, $H$-measurings $\triangleright_U\colon H\otimes U\to U$, $\triangleright_V\colon H\otimes V\to V$ and a $U$-$V$-bimodule $M = {}_U M_V$, a \emph{bimodule $H$-measuring} on $M$ relative to $\triangleright_U$ and $\triangleright_V$ is a $\genfd$-linear map
$\triangleright_M\colon H\otimes M\to M$
such that $h\triangleright_M (u.m.v) = (h_{(1)}\triangleright_U u).(h_{(2)}\triangleright_M m).(h_{(3)}\triangleright_V v)$ for all $u\in U$, $m\in M$ and $v\in V$. 

If $\triangleright_U$ and $\triangleright_V$ are compatible with 2-cocycles $\sigma$ and $\tau$ respectively, then we say that $\triangleright_M$ is \emph{$\sigma,\tau$-compatible} if for all $h,k\in H$ and $m\in M$
\begin{equation}
\label{eq:compbim}
h\triangleright_M (k\triangleright_M m) = \sum \sigma(h_{(1)},k_{(1)}).((h_{(2)}k_{(2)})\triangleright_M m).\tau^{-1}(h_{(3)},k_{(3)}).
\end{equation}
\end{defn}

\begin{prop}
For any two cleft right $H$-comodule algebras $(E,\rho_E,\gamma_E)$ and $(F,\rho_F,\gamma_F)$ and a right $H$-comodule $({}_E N_F,\rho_N)$ in the category of $E$-$F$-bimodules, the $\genfd$-submodule of $H$-coinvariants $M = N^{\mathrm{co}H}$ is an $E^{\mathrm{co}H}$-$F^{\mathrm{co}H}$-bimodule by restriction.

The $\genfd$-linear map
\begin{equation}
\triangleright_M\colon H\otimes M\to M,\quad h\otimes w\mapsto \gamma_E(h_{(1)}). w.\gamma^{-1}_F(h_{(2)})
\end{equation}
is a bimodule $H$-measuring relative to $\triangleright_{E^{\mathrm{co}H}}$ and $\triangleright_{F^{\mathrm{co}H}}$, and it is $\sigma_{\gamma_E},\sigma_{\gamma_F}$-compatible, where $\sigma_{\gamma_E}$ and $\sigma_{\gamma_F}$ are defined by formula~(\ref{eq:sigmagamma}).
\end{prop}

\begin{thm}
Suppose $(U,\triangleright_U,\sigma)$ and $(V,\triangleright_V,\tau)$ are crossed systems and $\triangleright_M$ a $\sigma,\tau$-compatible bimodule $H$-measuring on ${}_U M_V$ relative to $\triangleright_U$ and $\triangleright_V$.

(i) The $\genfd$-module $M\otimes H$ admits a structure of $U\sharp_\sigma H$-$V\sharp_\tau H$-bimodule,
denoted $M\sharp_{\sigma,\tau} H = M\sharp_{\triangleright_M,\sigma,\tau} H$,
with the left $U\sharp_\sigma H$-action
\begin{equation}\label{eq:leftact}
(u\sharp h).(m\otimes g) = u.(h_{(1)}\triangleright_M m).\tau(h_{(2)},g_{(1)})\otimes h_{(3)}g_{(2)}
\end{equation}
and the right $V\sharp_\tau H$-action
\begin{equation}\label{eq:rightact}
(m\otimes g).(v\sharp k) = m.\bigl((g_{(1)}\triangleright_V v)\tau(g_{(2)},k_{(1)})\bigr)\otimes g_{(3)}k_{(2)}.
\end{equation}

(ii) $M\sharp_{\sigma,\tau}H$ has the canonical right $H$-comodule structure $\id_M\otimes\Delta$ which makes it into an object in ${}_{U\sharp H}\mathcal{M}^H_{V\sharp H}$. That means that the right $H$-comodule structure $\id\otimes\Delta$ on $M\sharp_{\sigma,\tau}H$ is a map of $U\sharp_\sigma H$-$V\sharp_\tau H$-bimodules where the bimodule structure on $M\sharp H\otimes H$ is
\begin{equation}
  (u\sharp h).(m\otimes k\otimes k').(v\sharp l) = (u\sharp h_{(1)}).(m\otimes k).(v\sharp l_{(1)})\otimes h_{(2)} k'l_{(2)}.
\end{equation}
In other words, 
$$
(\id_M\otimes\Delta)((u\sharp h).(m\otimes k).(v\sharp l)) =
((\id_U\otimes\Delta)(u\sharp h)).((\id_M\otimes\Delta)(m\otimes k)).((\id_V\otimes\Delta)(v\sharp l)). 
$$
\end{thm}

\begin{proof}
(i) Bilinearity of the candidates given by formulas~(\ref{eq:leftact}) and (\ref{eq:rightact}) is obvious.
To check the left action axiom, we compute for all $u,\tilde{u}\in U$, $m\in M$ and $h, \tilde{h},g\in H$:
$$\begin{array}{l}
(u\sharp h).((\tilde{u}\sharp\tilde{h}).(m\otimes g)) 
     = (u\sharp h).(\tilde{u}.(\tilde{h}_{(1)}\triangleright_M m).\tau(\tilde{h}_{(2)},g_{(1)})\otimes \tilde{h}_{(3)}g_{(2)})\\
     = u.(h_{(1)}\triangleright_M (\tilde{u}.(\tilde{h}_{(1)}\triangleright_M m).\tau(\tilde{h}_{(2)},g_{(1)})).\tau(h_{(2)},\tilde{h}_{(3)}g_{(2)})\otimes h_{(3)}\tilde{h}_{(4)}g_{(3)}\\
     = u(h_{(1)}\triangleright_U\tilde{u}).(h_{(2)}\triangleright_M(\tilde{h}_{(1)}\triangleright_M m)).(h_{(3)}\triangleright_V\tau(\tilde{h}_{(2)},g_{(1)}))\tau(h_{(4)},\tilde{h}_{(3)}g_{(2)})\otimes h_{(5)}\tilde{h}_{(4)}g_{(3)}\\
     \stackrel{(\ref{eq:2coc})}=  u(h_{(1)}\triangleright_U\tilde{u}).(h_{(2)}\triangleright_M(\tilde{h}_1\triangleright_M m)).\tau(h_{(3)},\tilde{h}_{(2)})\tau(h_{(4)}\tilde{h}_{(3)},g_{(1)})\otimes h_{(5)}\tilde{h}_{(4)}g_{(2)}\\
     \stackrel{(\ref{eq:compbim})}= u(h_{(1)}\triangleright_U\tilde{u})\sigma(h_{(2)},\tilde{h}_{(1)}).((h_{(3)}\tilde{h}_{(2)})\triangleright_M m).\tau^{-1}(h_{(4)},\tilde{h}_{(3)})\tau(h_{(5)},\tilde{h}_{(4)})\tau(h_{(6)}\tilde{h}_{(5)},g_{(1)})\otimes h_{(7)}\tilde{h}_{(6)}g_{(2)}\\
     = u(h_{(1)}\triangleright_U\tilde{u})\sigma(h_{(2)},\tilde{h}_{(1)}).((h_{(3)}\tilde{h}_{(2)})\triangleright_M m).\tau(h_{(4)}\tilde{h}_{(3)},g_{(1)})\otimes h_{(5)}\tilde{h}_{(4)}g_{(2)} \\
     = (u(h_{(1)}\triangleright_U\tilde{u})\sigma(h_{(2)},\tilde{h}_{(1)})\sharp h_{(3)}\tilde{h}_{(2)}).(m\otimes g)\\
     = ((u\sharp h)(\tilde{u}\sharp\tilde{h})).(m\otimes g).
\end{array}$$
To check the right action axiom we proceed as
$$\begin{array}{l}
(m\otimes g).((v\sharp h)(w\sharp k)) 
     = (m\otimes g).(v(h_{(1)}\triangleright_V w)\tau(h_{(2)},k_{(1)})\sharp h_{(3)}k_{(2)})\\
     = m.(g_{(1)}\triangleright_V v)(g_{(2)}\triangleright_V(h_{(1)}\triangleright_V w))(g_{(3)}\triangleright_V \tau(h_{(2)},k_{(1)}))\tau(g_{(4)},h_{(3)}k_{(2)})\sharp g_{(5)}h_{(4)}k_{(3)}\\
     = m.(g_{(1)}\triangleright_V v)\tau(g_{(2)},h_{(1)})(g_{(3)}h_{(2)}\triangleright_V w)\tau^{-1}(g_{(4)},h_{(3)})\tau(g_{(5)},h_{(4)})\tau(g_{(6)}h_{(5)},k_{(1)})\sharp g_{(7)}h_{(6)}k_{(2)}\\
     = m.(g_{(1)}\triangleright_V v)\tau(g_{(2)},h_{(1)})(g_{(3)}h_{(2)}\triangleright_V w)\tau(g_{(4)}h_{(3)},k_{(1)})\sharp g_{(5)}h_{(4)}k_{(2)}\\
     = (m.(g_{(1)}\triangleright_V v)\tau(g_{(2)},h_{(1)})\otimes g_{(3)}h_{(2)}).(w\sharp k)\\
     = ((m\otimes g).(v\sharp h)).(w\sharp k).
\end{array}$$

The unitality of both the left and right actions is straightforward.

Finally, we check that the left and right actions commute:
$$\begin{array}{l}
(u\sharp h).((m\otimes g).(v\sharp k)) = (u\sharp h)(m.((g_{(1)}\triangleright_V v)\tau(g_{(2)},k_{(1)}))\otimes g_{(3)}k_{(2)}) \\
= u.(h_{(1)}\triangleright_M m).(h_{(2)}\triangleright_V(g_{(1)}\triangleright_V v))(h_{(3)}\triangleright_V\tau(g_{(2)},k_{(1)}))\tau(h_{(4)},g_{(3)}k_{(2)})\otimes h_{(5)}g_{(4)}k_{(3)}
\\
= u.(h_{(1)}\triangleright_M m).\tau(h_{(2)},g_{(1)})(h_{(3)}g_{(2)}\triangleright_V v)\tau^{-1}(h_{(4)},g_{(3)})\tau(h_{(5)},g_{(4)})\tau(h_{(6)}g_{(5)},k_{(1)})\otimes h_{(7)}g_{(6)}k_{(2)}
\\
= u.(h_{(1)}\triangleright_M m).\tau(h_{(2)},g_{(1)})(h_{(3)}g_{(2)}\triangleright_V v)\tau(h_{(4)}g_{(3)},k_{(1)})\otimes h_{(5)}g_{(4)}k_{(2)}
\\
=(u.(h_{(1)}\triangleright_M m).\tau(h_{(2)},g_{(1)})\otimes h_{(3)}g_{(2)}).(v\sharp k)\\
= ((u\sharp h).(m\otimes g)).(v\sharp k).
\end{array}$$
(ii) This is simple and left to the reader.
\end{proof}
\begin{prop}
  Let $(E,\rho^E,\gamma^E)$, $(F,\rho^F,\gamma^F)$ be cleft right $H$-comodule algebras and $M$ an object in ${}_E\mathcal{M}_F^H$. Consider the $H$-bimodule measuring $\triangleright = \triangleright_{M^{\mathrm{co}H}}$, $h\triangleright m = \gamma^E(h_{(1)})m(\gamma^F)^{-1}(h_{(2)})$. By restriction of scalars along isomorphisms of $H$-comodule algebras $E\cong E^{\mathrm{co}H}\sharp_{\triangleright_{E^{\mathrm{co}H}},\sigma_{\gamma^E}} H$ and $F\cong F^{\mathrm{co}H}\sharp_{\triangleright_{F^{\mathrm{co}H}},\sigma_{\gamma^F}} H$ the crossed product bimodule $M^{\mathrm{co}H}\sharp_{\triangleright,\sigma_{\gamma^E},\sigma_{\gamma^F}} H$ is viewed as an object in ${}_E\mathcal{M}_F^H$. Then
  \begin{equation}
    M^{\mathrm{co}H}\sharp_{\triangleright,\sigma_{\gamma^E},\sigma_{\gamma^F}} H\to M,\,\,\,\,
    m\sharp h\mapsto m\gamma^F(h)
  \end{equation}
  is an isomorphism in ${}_E\mathcal{M}_F^H$ with the inverse
  $m\mapsto m_{(0)}(\gamma^F)^{-1}(m_{(1)})\sharp m_{(2)}$.
\end{prop}

\subsection{Relation to $E\square^H F^{\mathrm{op}}$-modules}

If $H$ is a Hopf $\genfd$-algebra with a bijective antipode $S$, $(E,\rho^E)$ and $(F,\rho^F)$ (with $\rho_F(f) = f_{(0)}\otimes f_{(1)}$) right $H$-comodule algebras then make the opposite algebra $F^{\mathrm{op}}$ into a left $H$-comodule algebra via coaction $\rho' = (S^{-1}\otimes\id_F)\circ\tau\circ\rho^F$ where $\tau$ flips the tensor factors. Therefore the cotensor product algebra $E\square^H F^{\mathrm{op}}$ is defined using $\rho^E$ and $\rho'$. On the other hand, the algebra $E\otimes F^{\mathrm{op}}$ can be considered; it is an object of ${}_E\mathcal{M}^H_F$ via the right $H$-coaction $e\otimes f\mapsto e_{(0)}\otimes f_{(0)}\otimes e_{(1)}f_{(1)}$. There is a canonical isomorphism $(E\otimes F^{\mathrm{op}})^{\mathrm{co}H}=E\square^H F^{\mathrm{op}}$ (see~\cite{caenmarcustakeuchi}, Section 4).

In~\cite{caenmarcustakeuchi}
functor $(E\otimes F^{\mathrm{op}})\otimes_{E\square^H F^{\mathrm{op}}}-:{}_{E\square^H F^{\mathrm{op}}}\mathcal{M}\to{}_E\mathcal{M}^H_F$ and its right adjoint $(-)^{\mathrm{co}H}:{}_E\mathcal{M}^H_F\to{}_{E\square^H F^{\mathrm{op}}}\mathcal{M}$ have been exhibited. They prove that if $E$ is faithfully flat $H$-Hopf-Galois then this pair of functors forms an adjoint equivalence. In the setting of Subsection~\ref{ssec:DTbim}, if $(E,\rho^E,\gamma^E)$ and $F$ are cleft $H$-comodule algebras with fixed cleavings, then our $(-)^{\mathrm{co}H}$ takes value in a category $\operatorname{Cocyc}^H(E^{\mathrm{co}H},\triangleright_{\gamma^E},\sigma_{\gamma^E},F^{\mathrm{co}H},\triangleright_{\gamma^F},\sigma_{\gamma^F})$ of $E^{\mathrm{co}H}$-$F^{\mathrm{co}H}$-bimodules equipped with $\sigma_{\gamma^E},\sigma_{\gamma^F}$-compatible bimodule $H$-measurings relative to $\triangleright_{\gamma^E},\triangleright_{\gamma^F}$. Thus, if both $E$ and $F$ are cleft, one is lead to a correspondence between left $E\square^H F^{\mathrm{op}}$-modules and $E^{\mathrm{co}H}$-$F^{\mathrm{co}H}$-bimodules with compatible $H$-measurings, for a given choice of cleavings. Since $E^{\mathrm{co}H}\otimes(F^{\mathrm{co}H})^{\mathrm{op}}$ is a subalgebra of $E\square^H F^{\mathrm{op}}$, the restriction of scalars induces a faithful functor ${}_{E\square^H F^{\mathrm{op}}}\mathcal{M}\to{}_{E^{\mathrm{co}H}}\mathcal{M}_{F^{\mathrm{co}H}}$. The missing information is then imported into the bimodule measuring.
    \begin{lem}
If $E$ and $F$ are cleft then the $\genfd$-linear map $E^{\mathrm{co}H}\otimes_\genfd H\otimes_\genfd F^{\mathrm{co}H}\to E\otimes_\genfd F^{\mathrm{op}}$ extending the formula
$$u\otimes h\otimes v\mapsto u\gamma^E(h_{(1)})\otimes_\genfd v\cdot^{\mathrm{op}}(\gamma^F)^{-1}(h_{(2)}),\,\,\,\,\,u\in E^{\mathrm{co}H}, v\in F^{\mathrm{co}H}, h\in H$$ is an isomorphism onto its image which is precisely the subalgebra $E\square^H F^{\mathrm{op}}\hookrightarrow E\otimes_\genfd F^{\mathrm{op}}$. 
    \end{lem}
    \begin{proof}
  By an immediate calculation using $\rho_F((\gamma^F)^{-1}(k)) = (\gamma^F)^{-1}(k_{(2)})\otimes S(k_{(1)})$, each element of the form $u\gamma^E(h_{(1)})\otimes_\genfd v\cdot^{\mathrm{op}}(\gamma^F)^{-1}(h_{(2)})$ is coinvariant. The map $E^{\mathrm{co}H}\otimes H\otimes H\otimes F^{\mathrm{co}H}\to E\otimes F$ given by $u\otimes h\otimes k\otimes v\mapsto u\gamma^E(h)\otimes(\gamma^F)^{-1}(k)v$ is an isomorphism of $E^{\mathrm{co}H}$-$F^{\mathrm{co}H}$-bimodules (not algebras!) with an explicit inverse $e\otimes f\mapsto e_{(0)}(\gamma^{E})^{-1}(e_{(1)})\otimes e_{(2)}\otimes f_{(-2)}\otimes\gamma^F(f_{(-1)})f_{(0)}$. It is now sufficient to show that if $\sum_i\gamma^E(h_i)\otimes(\gamma^F)^{-1}(k_i)$ is a coinvariant in $E\otimes F$, then it is of the form $\gamma^E(h_{(1)})\otimes(\gamma^F)^{-1}(h_{(2)})$ for some $h\in H$; note that this is in the image $\gamma^E(H)\otimes(\gamma^F)^{-1}(H)$ of $H\otimes H = \genfd\otimes H\otimes H\otimes\genfd\subset E^{\mathrm{co}H}\otimes H\otimes H\otimes F^{\mathrm{co}H}$.
Applying the right $H$-coaction to the tensor product of $H$-comodules $\gamma^E(H)\otimes(\gamma^F)^{-1}(H)$ in terms of the corresponding elements in $H\otimes H$ we obtain $\sum_i h_i\otimes k_i\otimes 1 = \sum_i h_{i(1)}\otimes k_{i(2)}\otimes h_{i(2)}S(k_{i(1)})$ and therefore
      $\sum_i h_i\otimes k_i = \sum_i h_{i(1)}\otimes  h_{i(2)}S(k_{i(1)})k_{i(2)}
      = \sum_i (\epsilon(k_i)h_{i})_{(1)}\otimes  (\epsilon(k_i)h_{i})_{(2)}$.
    \end{proof}
    Thus, the additional information on $N\in {}_{E^{\mathrm{co}H}}\mathcal{M}_{F^{\mathrm{co}H}}$ is simply applying elements of the form $\gamma^E(h_{(1)})\otimes(\gamma^F)^{-1}(h_{(2)})$ on $N$ which amounts to the bimodule $H$-measuring $h\triangleright_N n = \gamma^E(h_{(1)}).n.(\gamma^F)^{-1}(h_{(2)})$.
\begin{thm}
  There is an isomorphism between category $\operatorname{Cocyc}^H(E^{\mathrm{co}H},\triangleright_{\gamma^E},\sigma_{\gamma^E},F^{\mathrm{co}H},\triangleright_{\gamma^F},\sigma_{\gamma^F})$ and category ${}_{E\square^H F^{\mathrm{op}}}\mathcal{M}$ which is identical on the underlying (restricted) $E^{\mathrm{co}H}$-$F^{\mathrm{co}H}$-bimodule $N$ with bimodule $H$-measuring $\triangleright_N$ and equips it with the left $E\square^H F^{\mathrm{op}}$-action
  $$
\left(\sum_i u_i\gamma^E(h_{i(1)})\square (v_i\cdot^{\mathrm{op}}(\gamma^{F})^{-1}(h_{i(2)}))\right)\otimes n\mapsto \sum_i u_i.(h_{i}\triangleright_N n).v_i
$$
for $u_i\in E^{\mathrm{co}H}$, $v_i\in F^{\mathrm{co}H}$, $h_i\in H$ and $n\in N$ where $\triangleright_N$ is the bimodule $H$-measuring on $N$ understood as an object in $\operatorname{Cocyc}^H(E^{\mathrm{co}H},\triangleright_{\gamma^E},\sigma_{\gamma^E},F^{\mathrm{co}H},\triangleright_{\gamma^F},\sigma_{\gamma^F})$.
\end{thm}
If we want to see the same picture in the case of crossed products, we recall that the inverse of the cleaving $\gamma_\sigma:h\mapsto 1\sharp_{\triangleright,\sigma} h = 1\sharp h$ is $\gamma_\sigma^{-1}:h\mapsto \sigma^{-1}(S h_{(2)},h_{(3)})\sharp S h_{(1)}$ hence
$$\begin{array}{lcl}
      (v\sharp 1)\cdot^{\mathrm{op}}\gamma^{-1}_\sigma(h) &=& (\sigma^{-1}(S h_{(2)},h_{(3)})\sharp S h_{(1)})\cdot (v\sharp 1) \\
      &=& \sigma^{-1}(S h_{(2)},h_{(3)}) ((S h_{(1)})_{(1)}\triangleright v)\sigma((S h_{(1)})_{(2)},1)\sharp (S h_{(1)})_{(3)}\\
  &=& \sigma^{-1}(S h_{(3)},h_{(4)}) (S h_{(2)}\triangleright v)\sharp S h_{(1)}.\end{array} $$
From this, we can explicitly present the cotensor product algebra $(U\sharp H)\square^H (V\sharp H)^{\mathrm{op}}$ in terms of $H$-measurings and 2-cocycles but without cotensor product, namely we send $u\gamma_\sigma(h_{(1)})\otimes v\cdot^{\mathrm{op}}\gamma_\sigma^{-1}(h_{(2)})$ to $u\otimes h\otimes v$ in the $\genfd$-module $M = U\otimes_\genfd H\otimes_\genfd V$ and transport the product. To prepare for this, we first observe that in any cleft $H$-comodule algebra $(F,\rho,\gamma)$, (thus, in particular, in a crossed product of the form $V\sharp_{\triangleright,\sigma}H$) the following formula holds for all $k,h\in H$ and $v',v\in F^{\mathrm{co}H}$:
\begin{equation}
  \gamma^{-1}(k)v'\gamma^{-1}(h)v =
  \gamma^{-1}(h_{(1)}k_{(1)})\sigma^{-1}_\gamma(h_{(2)},k_{(2)})(h_{(3)}\triangleright_{\gamma}v')v.
\end{equation}
    \begin{thm}
      For any $(U,\triangleright_U,\sigma^U)$ and $(V,\triangleright_V,\sigma_V)$,  the $\genfd$-module $M = U\otimes_\genfd H\otimes_\genfd V$ is an associative $\genfd$-algebra via the multiplication $\cdot_M$ on $M$ given by the bilinear extension of the formulas
        $$
      (u\otimes h\otimes v)\cdot_M(u'\otimes k\otimes v') =
      u (h_{(1)}\triangleright_U u')\sigma_U(h_{(2)},k_{(1)})\otimes h_{(3)}k_{(2)}
      \otimes \sigma^{-1}_V(h_{(4)},k_{(3)})(h_{(5)}\triangleright_V v') v
      $$

There is an isomorphism $(M,\cdot_M)\mapsto(U\sharp_{\triangleright_U,\sigma_U} H)\square^H(V\sharp_{\triangleright_V,\sigma_V} H)^{\mathrm{op}}$ of algebras given by the formula
      $$
u\otimes h\otimes v\mapsto  (u\sharp h_{(1)})\otimes((v\sharp 1)\cdot^{\mathrm{op}}\gamma^{-1}_\sigma(h_{(2)})) = (u\sharp h_{(1)})\otimes(\sigma^{-1}_V(S h_{(4)},h_{(5)})(Sh_{(3)}\triangleright_V v)\sharp Sh_{(2)})^{\mathrm{op}}.
$$
Here the opposite mark on the entire right hand tensor factor is just to emphasize that it is viewed as an element in the opposite algebra.

To each $(N,\triangleright_N)$ in $\operatorname{Cocyc}^H(U,\triangleright_U,\sigma^U,V,\triangleright_V,\sigma_V)$, one assigns a left $(M,\cdot_M)$-action
\begin{equation}\label{eq:actionUHV}
(u\otimes h\otimes v).n = u.(h\triangleright_N n).v
\end{equation}
This correspondence, with an identity on morphisms, becomes an isomorphism of
the category of left $(M,\cdot_M)$-modules and the category $\operatorname{Cocyc}^H(U,\triangleright_U,\sigma^U,V,\triangleright_V,\sigma^V)$.
    \end{thm}
    It may be useful for the reader to check the action axiom for the action~(\ref{eq:actionUHV}) directly, using~(\ref{eq:compbim}), the bimodule measuring axiom, and the formula for the product $\cdot_M$.
  
\section{Biequivalence of bicategories}

Fix a Hopf $\genfd$-algebra $H$. We define two bicategories which roughly look like the standard bicategory $\operatorname{Bimod}$ of algebras, bimodules and bimodule maps, with additional structure and conditions on the objects, morphisms and 2-cells involving $H$-coactions and cleavings in the first case, and $H$-measurings and 2-cocycles in the latter.

\subsection{Bicategories $\operatorname{Bimod}^H$ and $\operatorname{ClBimod}^H$}

The bicategory $\operatorname{Bimod}$ of rings, bimodules and morphisms of bimodules is internalized to $H$-comodules to get the bicategory $\operatorname{Bimod}^H$ of $H$-comodule algebras, $H$-comodule bimodules, and $H$-equivariant morphisms of bimodules. Furthermore, for the cleft case we introduce a refined bicategory $\operatorname{ClBimod}^H$.

Objects of $\operatorname{ClBimod}^H$ are cleft right $H$-comodule algebras $(E,\gamma)$.

The morphism category from $(E,\rho_E,\gamma_E)$ to $(F,\rho_F,\gamma_F)$ is the category ${}_E\mathcal{M}^H_F$; it does not depend on the cleavings $\gamma_E,\gamma_F$. An object in ${}_E\mathcal{M}^H_F$ is an $E$-$F$-bimodule $S$ with a right $H$-coaction $\rho_S\colon S\to S\otimes H$ such that $\rho_S(e.s.f) = \rho_E(e).\rho_S(s).\rho_F(f)$ for all $e\in E, s\in S, f\in F$, where the $E\otimes H$-$F\otimes H$-bimodule structure on the right hand side combines the actions on $S$ and the multiplication in $H$ factorwise; in Sweedler notation,
$\rho_S(e.s.f) = e_{(0)}.s_{(0)}.f_{(0)}\otimes e_{(1)}s_{(1)}f_{(1)}\in S\otimes H$.
Morphisms in ${}_E\mathcal{M}^H_F$ are maps of $E$-$F$-bimodules which are also maps of $H$-comodules. Descriptively, 1-cells in $\operatorname{ClBimod}^H$ are bimodules with a compatible right $H$-coaction and 2-cells are bimodule $H$-comodule maps, with vertical composition being composition of such maps.

Horizontal composition is the tensor product of bimodules, with the $H$-coaction $\genfd$-linearly extending formula 
$\rho_{S\otimes_F T}(s\otimes t) = (s_{(0)}\otimes t_{(0)})\otimes s_{(1)}t_{(1)}\in (S\otimes_F T)\otimes H$.

Coherences for horizontal composition are the coherences for the underlying tensor product of bimodules; they clearly respect the coaction.

\subsection{Bicategory of crossed systems $\operatorname{Cocyc}^H$}

Objects of $\operatorname{Cocyc}^H$ are crossed systems (Definition~\ref{def:crossedsystem}).

A morphism from $(U,\triangleright_U,\sigma_U)$ to $(V,\triangleright_V,\sigma_V)$ is a $U$-$V$-bimodule $P$ equipped with a bimodule $H$-measuring $\triangleright_P$ relative to $\triangleright_U,\triangleright_V$ and $\sigma_U,\sigma_V$-compatible. The identity morphism on $(U,\triangleright_U,\sigma_U)$ is $({}_U U_U,\triangleright_U)$.

A 2-cell from ${}_U P_V$ to ${}_U P'_V$ is a map $f\colon P\to P'$ of bimodules commuting with the bimodule $H$-measurings: $f(h\triangleright_{P} p) = h\triangleright_{P'} f(p)$ for all $h\in H$, $p\in P$.
Vertical composition of 2-cells is the composition of bimodule maps
(which automatically commutes with the measurings).

Horizontal composition is the tensor product of bimodules equipped with the tensor product of bimodule measurings.

\begin{defn}\label{def:tensmeas}
Given three algebras with $H$-measurings $(U,\triangleright_U)$, $(V,\triangleright_V)$, $(W,\triangleright_W)$, a bimodule $P = {}_U P_V$ with bimodule $H$-measuring $\triangleright_P$ relative to $\triangleright_U,\triangleright_V$, and a bimodule $Q = {}_V Q_W$ with bimodule $H$-measuring $\triangleright_Q$ relative to $\triangleright_V,\triangleright_W$, the \emph{tensor product $H$-measuring} is the $\genfd$-linear map
$$\triangleright_{P\otimes_V Q}\colon H\otimes (P\otimes_V Q)\to (P\otimes_V Q),
\quad h\otimes (p\otimes q)\mapsto (h_{(1)}\triangleright_P p)\otimes (h_{(2)}\triangleright_Q q).$$
\end{defn}

In the notation of Definition~\ref{def:tensmeas} the following holds:

\begin{prop} The tensor product measuring $\triangleright_{P\otimes_V Q}$ is a bimodule $H$-measuring relative to $\triangleright_U$ and $\triangleright_W$. If, in addition, $(U,\triangleright_U,\sigma_U)$, $(V,\triangleright_V,\sigma_V)$, $(W,\triangleright_W,\sigma_W)$ are crossed systems, $\triangleright_P$ is $\sigma_U,\sigma_V$-compatible and $\triangleright_Q$ is $\sigma_V,\sigma_W$-compatible, then $\triangleright_{P\otimes_V Q}$ is $\sigma_U,\sigma_W$-compatible.
\end{prop}

\subsection{Correspondence for tensor products of bimodules}

We now work with crossed systems $(U,\triangleright_U,\sigma)$, $(V,\triangleright_V,\tau)$, $(W,\triangleright_W,\upsilon)$ and $(Z,\triangleright_Z,\pi)$ and bimodules ${}_U P_V, {}_V Q_W, {}_W R_Z$ which are equipped with cocycle-compatible bimodule measurings. We create $P\sharp H$, $Q\sharp H$, $R\sharp H$ accordingly.
\begin{lem}\label{lem:tensoroversharp}
If $p\in {}_U P_V$, $q\in {}_V Q_W$, $v\in V$, $k,h,g\in H$, then in the $U$-$W$-bimodule $(P\sharp H)\otimes_{V\sharp_\tau H}(Q\sharp H)$
$$
\bigl( p.(g_{(1)}\triangleright_V v)\tau(g_{(2)},h_{(1)})\otimes g_{(3)}h_{(2)}\bigr)\otimes_{V\sharp_\tau H} (q\otimes k) = (p\otimes g)\otimes_{V\sharp_\tau H}
  \bigl(v.(h_{(1)}\triangleright_Q q).\upsilon(h_{(2)},k_{(1)})\otimes h_{(3)}k_{(2)}\bigr)
  $$
and in particular
  \begin{equation}\label{eq:tenoversharpsimplerform}
(p\otimes h)\otimes_{V\sharp_\tau H} (q\otimes k) = (p\otimes 1)\otimes_{V\sharp_\tau H}((h_{(1)}\triangleright_Q q).\upsilon(h_{(2)},k_{(1)})\otimes h_{(3)}k_{(2)})
  \end{equation}
\end{lem}

\begin{proof} The first identity follows immediately by applying~(\ref{eq:leftact}) and~(\ref{eq:rightact}) to $(p\otimes g).(v\sharp h)\otimes_{V\sharp_\tau H}(q\otimes k) = (p\otimes g)\otimes_{V\sharp_\tau H}(v\sharp h).(q\otimes k)$. Eq.~(\ref{eq:tenoversharpsimplerform}) is a special case when $v=1_V$ and $g = 1_H$.
\end{proof}

\begin{prop}\label{prop:cohforCross}
  There are canonical isomorphisms of $U\sharp H$-$W\sharp H$-bimodules and $H$-comodules $\zeta_{P,Q}:(P\otimes_V Q)\sharp H\cong (P\sharp H)\otimes_{V\sharp_\tau H}(Q\sharp H)$, $p\otimes q\otimes k\mapsto (p\otimes 1)\otimes_{V\sharp_\tau H}(q\otimes k)$ with well-defined inverse $(p\otimes h)\otimes_{V\sharp_\tau H}(q\otimes k)\mapsto (p\otimes_V (h_{(1)}\triangleright_Q q).\upsilon(h_{(2)},k_{(1)}))\otimes h_{(3)}k_{(2)}$.

  A coherence identity for $\zeta$ holds. Namely, the diagram
  \begin{equation}\label{eq:cohforCross}
\begin{tikzcd}[column sep=3em, row sep=5em]
  ((P\otimes_V Q)\otimes_W R)\sharp H
    \arrow[r, "\zeta_{P\otimes Q,\,R}"]
    \arrow[d, "\alpha_{P,Q,R}\,\sharp\, H"']
  & (P\otimes_V Q)\sharp H\otimes_{W\sharp H}R\sharp H
    \arrow[r, "\zeta_{P,Q}\otimes\mathrm{id}_{R\sharp H}"]
  & (P\sharp H\otimes_{V\sharp H}Q\sharp H)\otimes_{W\sharp H}R\sharp H
  \\
  (P\otimes_V (Q\otimes_W R))\sharp H
    \arrow[r, "\zeta_{P,\,Q\otimes R}"']
  & P\sharp H\otimes_{V\sharp H}(Q\otimes_W R)\sharp H
    \arrow[r, "\mathrm{id}_{P\sharp H}\otimes\,\zeta_{Q,R}"']
  & P\sharp H\otimes_{V\sharp H}
      \bigl(Q\sharp H\otimes_{W\sharp H}R\sharp H\bigr)
    \arrow[u, "\alpha_{P\sharp H,\,Q\sharp H,\,R\sharp H}"]
\end{tikzcd}
   \end{equation}
of $U\sharp H$-$Z\sharp H$-bimodules $H$-comodules commutes.
\end{prop}
\begin{proof}
  The fact that the formula for the inverse is well-defined, and that the two maps are indeed 2-sided inverses follows from the Lemma~\ref{lem:tensoroversharp}.

  The coherence identity is checked directly and easily on elements.
\end{proof}

Now we consider three cleft $H$-comodule algebras $E,F,G$ with cleavings $\gamma_E,\gamma_F,\gamma_G$, and bimodule $H$-comodules ${}_E M_F, {}_F N_G$. It is well-known that in that case the canonical map
\begin{equation}\label{eq:coinvtenprodbim}
M^{\mathrm{co}H}\otimes_{F^{\mathrm{co}H}}N^{\mathrm{co}H}\to(M\otimes_F N)^{\mathrm{co}H}
\end{equation}
is an isomorphism of $E^{\mathrm{co}H}$-$G^{\mathrm{co}H}$-bimodules. There are obvious coherence identities for three $\otimes$-composable bimodules.

\begin{prop}
The canonical isomorphism~(\ref{eq:coinvtenprodbim}) intertwines the corresponding bimodule $H$-measurings that is, it is a 2-cell in $\operatorname{Cocyc}^H$.
\end{prop}
\begin{proof}
  For any $h\in H$, $m\in M^{\mathrm{co}H}$, $n\in N^{\mathrm{co}H}$, the canonical isomorphism maps $(h_{(1)}\triangleright_M m)\otimes_{F^{\mathrm{co}H}}(h_{(2)}\triangleright_N n)\mapsto(h_{(1)}\triangleright_M m)\otimes_F(h_{(2)}\triangleright_N n) =
  \gamma_E(h_{(1)}).m.\gamma^{-1}_F(h_{(2)})\otimes_F \gamma_F(h_{(3)}).n.\gamma^{-1}_G(h_{(4)})=\gamma_E(h_{(1)}).m.\gamma^{-1}_F(h_{(2)})\gamma_F(h_{(3)})\otimes_F n.\gamma^{-1}_G(h_{(4)})=\gamma_E(h_{(1)}).m.\epsilon(h_{(2)})\otimes_F n.\gamma^{-1}_G(h_{(3)})=\gamma_E(h_{(1)}).m\otimes_F n.\gamma^{-1}_G(h_{(2)})= \gamma_E(h_{(1)}).(m\otimes_F n).\gamma^{-1}_G(h_{(2)}) = h\triangleright_{M\otimes_F N} (m\otimes_F n)$.
\end{proof}

\subsection{Doi--Takeuchi bicategory biequivalence}

We define a homomorphism of bicategories $\operatorname{Cross}^H\colon\operatorname{Cocyc}^H\to\operatorname{ClBimod}^H$:
on objects $(U,\triangleright,\sigma)\mapsto U\sharp_{\triangleright,\sigma} H$;
on morphisms $\big(({}_U P_V,\triangleright_P):(U,\triangleright_U,\sigma)\to(V,\triangleright_V,\tau)\big)\mapsto P\sharp_{\triangleright_P,\sigma,\tau} H$ with its canonical $H$-coaction;
on 2-cells $(f\colon {}_U P_V\to {}_U Q_V)\mapsto (f\otimes \id_H\colon P\sharp_{\triangleright_P,\sigma,\tau} H\to Q\sharp_{\triangleright_Q,\sigma,\tau} H)$.

\begin{remark}
  The coherence isomorphisms for compatibility of $\operatorname{Cross}^H$ with horizontal composition (tensor products) are the canonical isomorphisms of $U$-$V$-bimodules and $H$-comodules $(P\otimes_V Q)\otimes H\cong (P\otimes H)\otimes_{V\sharp_\tau H}(Q\otimes H)$ from Proposition~\ref{prop:cohforCross}.
\end{remark}

Homomorphism of bicategories $\operatorname{Coinv}^H\colon\operatorname{ClBimod}^H\to\operatorname{Cocyc}^H$ takes coinvariants and equips them with the induced measuring and cocycle data. More precisely, on objects $(E,\gamma)\mapsto(E^{\mathrm{co}H},\triangleright_\gamma,\sigma_\gamma)$ as in~\ref{ssec:crossedcleft}; on morphisms, any $H$-comodule $E$-$F$-bimodule $(S,\rho_S)\in{}_E\mathcal{M}^H_F$ is sent to $(S^{\mathrm{co}H},\triangleright_S)$, where $\triangleright_S$ is the bimodule $H$-measuring $h\otimes s\mapsto\gamma_E(h_{(1)}).s.\gamma_F^{-1}(h_{(2)})$ for $s\in S^{\mathrm{co}H}$; finally for every bimodule map $f:S\to S'$ the value $\operatorname{Coinv}^H(f) = f|_{S^{\mathrm{co}H}}:S^{\mathrm{co}H}\to S'^{\mathrm{co}H}$ is the restriction-corestriction of $f$ to the coinvariants. Noting that $\rho\circ\gamma_F^{-1}$ is the convolution inverse of $\rho\circ\gamma_F$, by its uniqueness, $\rho(\gamma_F^{-1}(h)) = \gamma_F^{-1}(h_{(2)})\otimes S h_{(1)}$ hence indeed $\rho(\gamma_E(h_{(1)}).s.\gamma_F^{-1}(h_{(2)})) = \gamma_E(h_{(1)}).s.\gamma_F^{-1}(h_{(4)})\otimes h_{(2)}S h_{(3)} = \gamma_E(h_{(1)}).s.\gamma_F^{-1}(h_{(2)})\otimes 1_H$ proving that the $\genfd$-submodule of $H$-coinvariants is a $E^{\mathrm{co}H}$-$F^{\mathrm{co}H}$-subbimodule. Measuring $\triangleright_S$ is $\sigma_{\gamma_E},\sigma_{\gamma_F}$-compatible by construction.

We now define pseudonatural 2-transformations
$\Phi\colon \mathrm{Id}_{\operatorname{Cocyc}^H}\Rightarrow \operatorname{Coinv}^H\circ \operatorname{Cross}^H$
and $\Psi\colon \mathrm{Id}_{\operatorname{ClBimod}^H}\Rightarrow \operatorname{Cross}^H\circ \operatorname{Coinv}^H$.

The component $\Phi_{(U,\triangleright,\sigma)}$ is the invertible bimodule ${}_U U_{U\otimes\genfd 1_H}={}_U U_{(U\sharp_{\triangleright,\sigma} H)^{\mathrm{co}H}}$ with right $U\otimes\genfd 1_H$-action $u.(w\otimes\lambda 1_H) = uw\lambda$ and the inverse bimodule $\Phi^{-1}_{(U,\triangleright,\sigma)}={}_{(U\sharp_{\triangleright,\sigma} H)^{\mathrm{co}H}}(U\sharp_{\triangleright,\sigma} H)^{\mathrm{co}H}_U$. It is viewed as an isomorphism of crossed systems from $(U,\triangleright,\sigma)$ to $((U\sharp_{\triangleright,\sigma}H)^{\mathrm{co}H},\triangleright_{\tilde\gamma},\sigma_{\tilde\gamma})$, where $U\cong U\otimes\genfd 1_H=(U\sharp H)^{\mathrm{co}H}\subset U\sharp_{\triangleright,\sigma} H$ and $\tilde\gamma(h) = 1\sharp h$. For a morphism $({}_U P_V,\triangleright_P)\colon(U,\triangleright_U,\sigma)\to(V,\triangleright_V,\tau)$ in $\operatorname{Cocyc}^H$, the component $\Phi_P = \Phi_{({}_U P_V,\triangleright_P)}$ fills the commutative square
$$
\begin{tikzcd}
(U,\triangleright_U,\sigma) \ar[r,"U"] \ar[d,"P"']\ar[dr,phantom,"\stackrel{\Phi_P}\Rightarrow"] &  (U\otimes\genfd 1_H,\triangleright_{\tilde\gamma},\sigma_{\tilde\gamma}) \ar[d,"P\otimes\genfd 1_H"]\\
(V,\triangleright_V,\tau)\ar[r,"V"'] & (V\otimes\genfd 1_H,\triangleright_{\tilde\gamma'},\tau_{\tilde\gamma'}),
\end{tikzcd}
$$

\noindent namely, it is the canonical isomorphism of $U$-$V$-bimodules $P\otimes_V V\xrightarrow{\;\sim\;} U\otimes_{U\otimes\genfd 1_H}(P\otimes_\genfd\genfd 1_H)$ given by (the $\genfd$-linear extension of) $p \otimes v \mapsto 1_U \otimes (p.v \otimes 1_H)$ with inverse $u\otimes p\otimes 1_H \mapsto u.p \otimes 1_V$.

\begin{lem}
  $\Phi\colon \mathrm{Id}_{\operatorname{Cocyc}^H}\Rightarrow \operatorname{Coinv}^H\circ \operatorname{Cross}^H$ satisfies the coherences for a pseudonatural transformation.
\end{lem}
\begin{proof} The naturality in 2-cells means that for any 2-cell $f:({}_U P_V,\triangleright,\sigma_U,\sigma_V)\to({}_U P'_V,\triangleright,\sigma'_U,\sigma'_V)$ in $\operatorname{Cocyc}^H$ the identity
$(f\sharp\id_H)|_{(P\sharp H)^{\mathrm{co}H}}\circ\Phi_P = \Phi_{P'}\circ (f\otimes\id_V):P\otimes_V V\to U\otimes_{U\otimes\genfd 1_H}(P'\otimes\genfd 1_H)$ holds. Both sides evaluate to $p\otimes v\mapsto 1_U\otimes (f(p.v)\sharp 1_H)$.

  We now show the compatibility of $\Phi$ with the horizontal composition. The
  compatibility with the identity 2-cells is left to the reader.

Let $P \in \operatorname{Cocyc}^H(U, V)$ and $Q \in \operatorname{Cocyc}^H(V, W)$. The coherence for horizontal composition requires the following diagram to commute in the category of $U$-$W \otimes \genfd 1_H$ bimodules:

\begin{equation}\label{eq:PhiCoh}
\begin{tikzcd}[column sep=6em, row sep=5em]
(P\otimes_V Q) \otimes_W W \arrow[r, "\Phi_{P\otimes_V Q}"] \arrow[d, "\alpha"'] 
  & U \otimes_{U\otimes\genfd 1_H} ((P\otimes_V Q)\otimes\genfd 1_H) \arrow[d, "\cong"']\\
P \otimes_V (Q \otimes_W W) \arrow[d, "\id_P \otimes \Phi_Q"'] 
  & U \otimes_{U\otimes\genfd 1_H} ((P\otimes\genfd 1_H) \otimes_{V\otimes\genfd 1_H} (Q\otimes\genfd 1_H))  \\
P \otimes_V (V \otimes_{V\otimes\genfd 1_H} (Q\otimes\genfd 1_H)) \arrow[r, "\alpha^{-1}"'] 
  & (P \otimes_V V) \otimes_{V\otimes\genfd 1_H} (Q\otimes\genfd 1_H) \arrow[u, "\Phi_P \otimes \id"']
\end{tikzcd}
\end{equation}

To show the commutativity we trace an element $(p \otimes q) \otimes w \in (P\otimes_V Q) \otimes_W W$ through both paths from $(P\otimes_V Q) \otimes_W W$ to
$U \otimes_{U\otimes\genfd 1_H} ((P\otimes\genfd 1_H) \otimes_{V\otimes\genfd 1_H} (Q\otimes\genfd 1_H))$. The top arrow to the right applies $\Phi_{P\otimes_V Q}$, thus we multiply the elements in the tensor product over $W$ to obtain
$1_U \otimes ((p \otimes q.w) \otimes 1_H)$. The upper arrow of the right-hand column is induced by the standard tensor product identification $(p\otimes q) \otimes 1_H\mapsto (p\otimes 1_H) \otimes_{V\otimes\genfd 1_H} (q\otimes 1_H)$. Factoring our element through this gives $1_U \otimes ((p \otimes 1_H) \otimes_{V\otimes\genfd 1_H} (q.w \otimes 1_H))$.

The other path first applies associativity isomorphism $(p \otimes q) \otimes w \longmapsto p \otimes (q \otimes w)$, then $\id_P \otimes \Phi_Q$ to obtain
$p \otimes (1_V \otimes (q.w \otimes 1_H))$, again associativity
to get $(p \otimes 1_V) \otimes_{V\otimes\genfd 1_H} (q.w \otimes 1_H)$. 
Finally, apply $\Phi_P \otimes \id$ moving up the right edge
to obtain $(1_U \otimes (p.1_V \otimes 1_H)) \otimes_{V\otimes\genfd 1_H} (q.w \otimes 1_H)$. Because $p.1_V = p$, this simplifies to
$1_U \otimes ((p \otimes 1_H) \otimes_{V\otimes\genfd 1_H} (q.w \otimes 1_H))$.

Thus the diagram is commutative.
\end{proof}

The component $\Psi_{(E,\rho,\gamma)}\colon(E,\rho,\gamma)\to(E^{\mathrm{co}H}\sharp_{\triangleright_\gamma,\sigma_\gamma}H,\id\otimes\Delta,\tilde\gamma)$ is the invertible bimodule ${}_E E_{E^{\mathrm{co}H}\sharp_{\triangleright_\gamma,\sigma_\gamma}H}$ where the right action is $e.(u\sharp h) = eu\gamma(h)$. This is indeed a right action:
$$\begin{array}{lcl}
  (e.(u\sharp h)).(v\sharp k) &=& eu\gamma(h)v\gamma(k) = eu\gamma(h_{(1)})v\gamma^{-1}(h_{(2)})\gamma(h_{(3)})\gamma(k) \\
  &=&
  eu(h_{(1)}\triangleright_\gamma v)\sigma_\gamma(h_{(2)},k_{(1)})\gamma(h_{(3)}k_{(2)}) \\
  &=& e.(u(h_{(1)}\triangleright_\gamma v)\sigma_\gamma(h_{(2)},k_{(1)})\sharp h_{(3)}k_{(2)}) = e.((u\sharp h)(v\sharp k)).
  \end{array}$$
  The inverse bimodule is ${}_{E^{\mathrm{co}H}\sharp_{\triangleright_\gamma,\sigma_\gamma}H}(E^{\mathrm{co}H}\sharp_{\triangleright_\gamma,\sigma_\gamma}H)_E$ with the left regular action and right action given by $(u\sharp h).u'\gamma(h') = (u\sharp h)\cdot(u'\sharp h')$.
For a morphism ${}_E S_F : (E,\rho,\gamma_E)\to(F,\rho',\gamma_F)$,
the component $\Psi_S$ is the canonical isomorphism 
$$\begin{tikzcd}
  (E,\rho,\gamma_E) \ar[r,"E"] \ar[d,"S"']\ar[dr,phantom,"\stackrel{\Psi_S}\Rightarrow"] &  (E^{\mathrm{co}H}\sharp_{\triangleright_{\gamma_E},\sigma_{\gamma_E}}H,\id\otimes\Delta,\tilde\gamma_E) \ar[d,"S\otimes\genfd 1_H"]\\
(F,\rho',\gamma_F)\ar[r,"F"'] & (F^{\mathrm{co}H}\sharp_{\triangleright_{\gamma_F},\sigma_{\gamma_F}}H,\id\otimes\Delta,\tilde{\gamma}_F).
\end{tikzcd}$$

\begin{lem}
$\Psi\colon \mathrm{Id}_{\operatorname{ClBimod}^H}\Rightarrow \operatorname{Cross}^H\circ \operatorname{Coinv}^H$ is a pseudonatural transformation.
\end{lem}
\begin{proof}

For an object $(E, \rho_E, \gamma_E)$, the component $\Psi_E$ is the bimodule ${}_E E_{E^{\mathrm{co}H}\sharp H}$ where the right action is $e.(u\sharp h) = eu\gamma(h)$.
For a morphism $S = {}_E S_F \colon (E,\rho_E,\gamma_E)\to(F,\rho_F,\gamma_F)$, the canonical isomorphism $\Psi_S$ is explicitly

$$\Psi_S \colon S \otimes_F F \xrightarrow{\cong} E \otimes_{E^{\mathrm{co}H}\sharp H} (S^{\mathrm{co}H}\sharp H)$$

$$s \otimes f \longmapsto 1_E \otimes (s_{(0)}\gamma_F^{-1}(s_{(1)}) \sharp s_{(2)}) \cdot f$$

Because $S \otimes_F F \cong S$, this essentially maps $s \mapsto 1_E \otimes (s_{(0)}\gamma_F^{-1}(s_{(1)}) \sharp s_{(2)})$.

Let ${}_E M_F \colon E \to F$ and ${}_F N_G \colon F \to G$ be 1-cells in $\operatorname{ClBimod}^H$. The pseudonaturality coherence for horizontal composition requires the following diagram to commute, utilizing the coherence isomorphism $\zeta$ from Proposition~\ref{prop:cohforCross}:

\begin{equation}\label{eq:PsiCoh}
\begin{tikzcd}[column sep=7em, row sep=6em]
(M\otimes_F N) \otimes_G G \arrow[r, "\Psi_{M\otimes_F N}"] \arrow[d, "\alpha"'] 
  & E \otimes_{E^{\mathrm{co}H}\sharp H} ((M\otimes_F N)^{\mathrm{co}H} \sharp H)
  \arrow[d, "\id_E \otimes \zeta_{M,N}"']\\
M \otimes_F (N \otimes_G G) \arrow[d, "\id_M \otimes \Psi_N"'] 
  & E \otimes_{E^{\mathrm{co}H}\sharp H} ((M^{\mathrm{co}H} \sharp H) \otimes_{F^{\mathrm{co}H}\sharp H} (N^{\mathrm{co}H} \sharp H))  \\
M \otimes_F (F \otimes_{F^{\mathrm{co}H}\sharp H} (N^{\mathrm{co}H}\sharp H)) \arrow[r, "\alpha^{-1}"'] 
  & (M \otimes_F F) \otimes_{F^{\mathrm{co}H}\sharp H} (N^{\mathrm{co}H}\sharp H) \arrow[u, "\Psi_M \otimes \id"']
\end{tikzcd}
\end{equation}

Using the cleaving maps $\gamma_F$ and $\gamma_G$, any elements $s\in M$ and $t \in N$ can be uniquely written as $\genfd$-linear combinations $s = \sum_i m_i \cdot \gamma_F(h_i)$ and $t = \sum_j n_j\cdot \gamma_G(k_j)$ for some $m_i \in M^{\mathrm{co}H}, n_j\in N^{\mathrm{co}H}$ and $h_i, k_j \in H$.
Let us trace the element $(m_i \gamma_F(h_i) \otimes n_j \gamma_G(k_j)) \otimes 1_G$ along both paths. If we start to the bottom then first associativity gives
$(m_i\gamma_F(h_i) \otimes_F (n_j \gamma_G(k_j) \otimes 1_G)$.
Apply $\id_M \otimes \Psi_N$. Observing that $n$ is coinvariant, $\Psi_N(n_j \gamma_G(k_j) \otimes 1_G) = 1_F \otimes (n_j \sharp k_j)$, hence we obtain
$m_i \gamma_F(h_i) \otimes_F (1_F \otimes_{F^{\mathrm{co}H}\sharp H} (n_j \sharp k_j))$,
and after associativity $(m_i \gamma_F(h_i) \otimes_F 1_F) \otimes_{F^{\mathrm{co}H}\sharp H} (n_j \sharp k_j)$. Since $m_i \gamma_F(h_i) \otimes_F 1_F$ is identified with $m_i \otimes_F \gamma_F(h_i)$ (moving the scalar across the tensor product over $F$), we apply $\Psi_M \otimes \id$ upward. $\Psi_M$ acting on $m_i\otimes_F\gamma_F(h_i)$ evaluates to the right action of $\gamma_F(h_i)$ on $1_E \otimes (m_i\sharp 1_H)$, which in the crossed product $(M^{\mathrm{co}H}\sharp H)$ becomes $1_E \otimes (m_i\sharp h_i)$. Thus, the right edge yields:
$1_E \otimes ((m_i\sharp h_i) \otimes_{F^{\mathrm{co}H}\sharp H} (n_j\sharp k_j))$.

For computing the top arrow to the right, we first rewrite
$(m_i\gamma_F(h_i)\otimes_F n_j\gamma_G(k_j)) \otimes 1_G = m_i\otimes_F \gamma_F(h_i)n_j\gamma_G(k_j)
= m_i\otimes_F (h_{i(1)} \triangleright_N n_j) \gamma_G(h_{i(2)}) \gamma_G(k_j)
= \left[ m \otimes_F (h_{(1)} \triangleright_N n) \sigma_G(h_{(2)}, k_{(1)}) \right] \cdot \gamma_G(h_{(3)}k_{(2)})$. The part in the square bracket is in $(M\otimes_F N)^{\mathrm{co}H}$. Thus applying $\Psi_{M\otimes_F N}$ results in 
$1_E\otimes\left(m_i\otimes_F (h_{i(1)}\triangleright_N n_j) \sigma_G(h_{i(2)}, k_{j(1)}) \sharp h_{i(3)} k_{j(2)} \right)$. Then apply $\id_E\otimes\zeta^{-1}_{M,N}$, where the formula for the inverse of the map $\zeta$ is from the statement of Proposition~\ref{prop:cohforCross}). This yields $1_E \otimes ((m_i\sharp h_i) \otimes_{F^{\mathrm{co}H}\sharp H} (n_j\sharp k_j))$. Thus both paths evaluate to the same element.
\end{proof}
\begin{thm}
  Homomorphisms of bicategories
  $$\operatorname{Cross}^H\colon\operatorname{Cocyc}^H\to\operatorname{ClBimod}^H\,\,\,\,\mathrm{ and }\,\,\,\,\operatorname{Coinv}^H\colon\operatorname{ClBimod}^H\to\operatorname{Cocyc}^H,$$
together with pseudonatural equivalences
$$\Phi\colon \mathrm{Id}_{\operatorname{Cocyc}^H}\Rightarrow \operatorname{Coinv}^H\circ \operatorname{Cross}^H
\,\,\,\,\mathrm{ and }\,\,\,\,
\Psi\colon \mathrm{Id}_{\operatorname{ClBimod}^H}\Rightarrow \operatorname{Cross}^H\circ \operatorname{Coinv}^H$$
form a biequivalence of bicategories. 
\end{thm}

\subsection{Correspondences for corings and noncommutative covers}

Comonads in the standard bicategory $\operatorname{Bimod}$ of rings, bimodules and bimodule maps are known as corings~\cite{brzwisCorings}. In any monoidal category admitting coequalizers of parallel pairs which commute with the tensor product, one has relative tensor products over monoids and the bicategory of monoids, bimodules and bimodule maps is defined; comonads in this category are known as internal corings. Since the bicategories $\operatorname{Bimod}^H$, $\operatorname{ClBimod}^H$ and $\operatorname{Cocyc}^H$ are equipped with faithful forgetful functors to $\operatorname{Bimod}$, comonads there are corings with additional structure. 

\subsubsection{Internal corings in $H$-comodules and $H$-Hopf modules over them}

Comonads in $\operatorname{Bimod}^H$ may also be thought of as internal corings in the monoidal category of $H$-comodules. Such an internal coring $D$ is described by a datum $(E,\rho^E,D,\rho^D,\Delta^D,\epsilon^D)$ which is simply a usual $E$-coring such that $(E,\rho^E)$ is a right $H$-comodule algebra, $(D,\rho^D)$ an object of the category ${}_E\mathcal{M}^H_E$ ($H$-comodules $E$-bimodule), and the maps $\Delta^D$ and $\epsilon^D$ morphisms of right $H$-comodules.

Recall that a (relative) left-right $(E,H)$-Hopf module $M = ({}_E M,\rho^M)$ is a left $E$-module ${}_E M$ with a right $H$-coaction $\rho^M$ which is a morphism of left $E$-modules, $\rho^M(em) = \rho^E(e)\rho^M(m)$. The category of $(E,H)$-Hopf modules is usually denoted by ${}_E\mathcal{M}^H$. For an $E$-coring $D$, a {\bf left-right $(D,H)$-Hopf comodule} is a triple $({}_E M,\delta^M,\rho^M)$ where
\begin{itemize}
\item $({}_E M,\rho^M)$ is a left-right $(E,H)$-Hopf module
\item $({}_E M, \delta^M)$ is a left $D$-comodule
\item the left $D$-coaction $\delta^M$ is a morphism in the category of $H$-comodules.
\end{itemize}
A morphism of left-right $(D,H)$-comodules is a morphism of left $D$-comodules which is also a morphism of right $H$-comodules. Left-right $(D,H)$-Hopf comodules and their morphisms form a category ${}^D\mathcal{M}^H$.

There is a comonad ${}_EG:{}_E\mathcal{M}\to{}_E\mathcal{M}$ which sends ${}_E M\mapsto {}_E (M\otimes H)$ where the left $E$-module structure on $M\otimes H$ is given by $e\otimes (m\otimes h)\mapsto e_{(0)}m\otimes e_{(1)}h$ and the component of the coproduct 
is $\id_M\otimes\Delta^H$. Comonad ${}_E G$ may be considered as the lift of the usual tensoring with $H$-comonad on the category of $\genfd$-modules $\mathcal{M}$ using the distributive law between $E\otimes -$ monad and $-\otimes H$ comonad on $\mathcal{M}$ (or, equivalenty, the entwining structure ${}_E\psi:E\otimes H\to H\otimes E$). It is well-known that the category of $(E,H)$-Hopf modules is canonically equivalent to the category of ${}_E G$-comodules. 

\begin{prop} There is a lift ${}^D G$ of the comonad ${}_E G$ along the forgetful functor ${}^D\mathcal{M}^H\to{}_E\mathcal{M}^H$ which is induced by the distributive law $\lambda : D\otimes_E({}_E G(-))\Rightarrow {}_E G(D\otimes_E -)$ of comonad $D\otimes_E-$ over comonad ${}_E G$ with components $\lambda_M:D\otimes {}_E(M\otimes H)\to {}_E(D\otimes_E M)\otimes H$ in ${}_E\mathcal{M}$ given by
\begin{equation}
\lambda_M(d\otimes_E(m\otimes h)) = (d_{(0)}\otimes_E m)\otimes d_{(1)}h
\end{equation}
This means that for a left $D$-comodule $M$ with $D$-coaction $m\mapsto m^{[-1]}\otimes m^{[0]}$, ${}^D G(M)$ is the left $E$-module ${}_E(M\otimes H)$ equipped with $D$-coaction $m\otimes h\mapsto \lambda_M((m^{[-1]}\otimes_E m^{[0]})\otimes h) = m^{[-1]}_{(0)}\otimes_E (m^{[0]}\otimes m^{[-1]}_{(1)}h)$.
\end{prop}

\begin{proof}
  First we check that the formula for $\lambda_M$ is well-defined on relative tensor product $\otimes_E$. Indeed, $\lambda_M(d.e \otimes m\otimes h) = ((d.e)_{(0)}\otimes_E m)\otimes (d.e)_{(1)} h = (d_{(0)}\otimes_E e_{(0)}m)\otimes d_{(1)}e_{(1)}h = \lambda_M(d\otimes e_{(0)} m\otimes e_{(1)}h)$. For the compatibility with the coproduct of $D$,
  $$(\Delta^D\otimes_E G(M))\circ\lambda_M =  \lambda_{D\otimes_E G(M)}\circ G(\lambda_{D\otimes_E M})\circ G(\Delta^D\otimes_E M):(D\otimes_E M)\otimes H\to D\otimes_E D\otimes_E (M\otimes H),$$
we check both sides on $(d\otimes_E m)\otimes h$. The left-hand side is $(d_{(0)})^{[1]}\otimes_E (d_{(0)})^{[2]}\otimes_E (m\otimes d_{(1)} h)$ and the right-hand side is $(d^{[1]})_{(0)}\otimes_E (d^{[2]})_{(0)}\otimes_E m\otimes (d^{[1]})_{(1)}(d^{[2]})_{(1)}h$ hence this follows by recalling that $\Delta^D:D\to D\otimes_E D$ is also a morphism in the monoidal category of $H$-comodules, hence $(d_{(0)})^{[1]}\otimes_E (d_{(0)})^{[2]}\otimes d_{(1)}  = (d^{[1]})_{(0)}\otimes_E (d^{[2]})_{(0)}\otimes (d^{[1]})_{(1)}(d^{[2]})_{(1)}$.
  Compatibility of $\lambda$ with $\delta^G$ and
  with counits is left to the reader.
\end{proof}

\begin{prop}
There is a canonical isomorphism between the category of left-right $(D,H)$-comodules and the category of ${}^D G$-comodules which is identical on morphisms and commuting with the forgetul functors to $E$-modules. On the level of objects, for any underlying left $E$-module $M$, left $D$-coactions and right $H$-coactions on $M$ making it into an object in ${}^D\mathcal{M}^H$ are in a bijective correspondence with ${}_E G$-coactions on $M$ which are simultaneously ${}^D G$-coactions.
\end{prop}

If $E$ is faithfully flat Hopf-Galois extension, the category of left-right $(D,H)$-Hopf comodules is equivalent to the category of $D^{\mathrm{co}H}$-comodules~\cite{skoda:hgcoring}, where ordinary coring $D^{\mathrm{co}H}$ is defined in~\ref{ssec:DTcmd}. This noncommutative descent along torsors theorem is a coring generalization of Schneider's theorem~\cite{Schn:prinHomSp} of the descent of $(E,H)$-Hopf modules for usual Hopf-Galois extensions and admits further generalizations.

\subsubsection{Comonads in $\operatorname{ClBimod}^H$}

A comonad in $\operatorname{ClBimod}^H$ consists of the following datum $(E,\rho^E,\gamma,D,\rho^D,\Delta^D,\epsilon^D)$:
\begin{itemize}
\item
    $E = (E,\rho^E)$ is a right $H$-comodule algebra with cleaving $\gamma\colon H\to E$
\item
    $D$ is an $E$-bimodule with right $H$-coaction $\rho^D$ making it an object of ${}_E\mathcal{M}_E^H$
\item
  $\Delta^D\colon D\to D\otimes_E D$, $d\mapsto d_{(1)}\otimes_E d_{(2)}$,
  is a map of $E$-bimodules and right $H$-comodules, where the $H$-coaction
  on $D\otimes_E D$ is the tensor product coaction
  $d\otimes_E d'\mapsto (d_{(0)}\otimes_E d'_{(0)})\otimes d_{(1)}d'_{(1)}$;
  moreover $\Delta^D$ should be coassociative,
  $(\id_D\otimes_E\Delta^D)\circ\Delta^D = (\Delta^D\otimes_E\id_D)\circ\Delta^D$, and counital via a counit $\epsilon^D\colon(D,\rho^D)\to (E,\rho^E)$,
  where $\epsilon^D$ is a morphism in ${}_E\mathcal{M}_E^H$.
\end{itemize}

\subsubsection{Comonads in $\operatorname{Cocyc}^H$}

A comonad in $\operatorname{Cocyc}^H$ is described by a datum
$(U,\triangleright_U,\sigma,C,\triangleright_C,\Delta^C,\epsilon^C)$ where:
\begin{itemize}
\item $(U,\triangleright_U,\sigma)$ is a crossed system
\item $C$ is a $U$-bimodule
\item $\triangleright_C\colon H\otimes C\to C$ is a bimodule $H$-measuring relative to $\triangleright_U,\triangleright_U$ and $\sigma,\sigma$-compatible
\item $(C,\Delta^C\colon C\to C\otimes_U C,\epsilon^C\colon C\to U)$ is a $U$-coring
\item $\Delta^C$ and $\epsilon^C$ are 2-cells in $\operatorname{Cocyc}^H$, i.e.\ they respect the measurings
\end{itemize}

\subsubsection{Doi--Takeuchi correspondence for comonads}
\label{ssec:DTcmd}

Any biequivalence of bicategories induces a biequivalence of the corresponding bicategories of comonads. Concretely, the bicategorical Doi--Takeuchi correspondence gives a biequivalence between the bicategories of comonads in $\operatorname{ClBimod}^H$ and comonads in $\operatorname{Cocyc}^H$; the homomorphisms of bicategories induced by $\operatorname{Cross}^H$ and $\operatorname{Coinv}^H$ on the bicategories of comonads will be denoted by $\operatorname{Cross}^H_{\mathrm{cmd}}$ and $\operatorname{Coinv}^H_{\mathrm{cmd}}$. 

For brevity we describe the object part of this biequivalence.
$\operatorname{Coinv}^H_{\mathrm{cmd}}(E,\rho^E,\gamma,D,\rho^D,\Delta^D,\epsilon^D)=(E^{\mathrm{co}H},\triangleright_{\gamma},\sigma_\gamma,D^{\mathrm{co}H},\Delta^{D^{\mathrm{co}H}},\epsilon^{D^{\mathrm{co}H}})$. The nontrivial part is the comultiplication $\Delta^{\mathrm{co}H}$ which may be defined in a satisfactory way for faithfully flat Hopf-Galois extensions.

  \begin{thm}\label{thm:coinvcor}
    Let $H$ be a Hopf algebra, and $(E,\rho^E)$ a faithfully flat right $H$-Hopf-Galois extension with subalgebra of $\rho^E$-coinvariants $E^{\mathrm{co}H}$ and $(D,\rho^D,\Delta^D,\epsilon^D)$ an $E$-coring in the monoidal category of right $H$-comodules.

    Then $D^{\mathrm{co}H}$ is an $E^{\mathrm{co}H}$-subbimodule of $D$ and the composition
    \begin{equation}
D^{\mathrm{co}H}\hookrightarrow D\to D\otimes_E D
    \end{equation}
    is a map of $E^{\mathrm{co}H}$-bimodules which corestricts to
    a map of $E^{\mathrm{co}H}$-bimodules $\Delta_1:D^{\mathrm{co}H}\hookrightarrow D\to (D\otimes_E D)^{\mathrm{co}H}$. The canonical map 
    $$
\nu_1 : D^{\mathrm{co}H}\otimes_{E^{\mathrm{co}H}}D^{\mathrm{co}H}\longrightarrow (D\otimes_E D)^{\mathrm{co}H}
    $$
is an isomorphism of $E^{\mathrm{co}H}$-bimodules.

The composition 
    $$
\Delta^{D^{\mathrm{co}H}}:D^{\mathrm{co}H}\stackrel{\Delta_1}\longrightarrow (D\otimes_E D)^{\mathrm{co}H}\stackrel{\nu_1^{-1}}\longrightarrow D^{\mathrm{co}H}\otimes_{E^{\mathrm{co}H}} D^{\mathrm{co}H}
  $$
is a coassociative comultiplication in the monoidal category of $E^{\mathrm{co}H}$-bimodules with a counit $\epsilon^{D^{\mathrm{co}H}}:D^{\mathrm{co}H}\to E^{\mathrm{co}H}$ given by the restriction-corestriction of $\epsilon^D:D\to E$.
  \end{thm}
The quasiinverse $\operatorname{Cross}^H_{\mathrm{cmd}}$ maps $(U,\triangleright_U,\sigma,C,\Delta^C,\epsilon^C)$ to the internal $U\sharp_{\triangleright_U,\sigma}H$-coring $C\sharp_{\triangleright,\sigma,\sigma}H$ with comultiplication $\Delta^C\sharp\id_H$ and $H$-coaction $\id_C\otimes\Delta$.

\subsubsection{Noncommutative covers}

One can additionally equip comonads in the bicategories $\operatorname{Bimod}$, $\operatorname{Bimod}^H$ and $\operatorname{ClBimod}^H$ with a map of comonads from the Sweedler comonad and require that this map is an epimorphism of underlying modules. For a $\genfd$-algebra $U$, the Sweedler coring corresponding to the extension $\genfd\hookrightarrow U$ is $U\otimes U$ with the external action $u.(v\otimes w).u' = uv\otimes wu'$ and comultiplication $u\otimes v\mapsto u\otimes 1 \otimes_U 1\otimes v\in (U\otimes U)\otimes_U(U\otimes U)$; in particular $1\otimes 1\in U\otimes U$ is grouplike. If $C$ is a $U$-coring, a surjective morphism $s:U\otimes U\to C$ of $U$-corings is determined by the choice of a grouplike element $g$ by setting $s(u\otimes v) = u.g.v$ or equivalently $g = s(1\otimes 1)$. Now, for an internal $E$-coring $D$ in the category of right $H$-comodules, $s=s^D$ is in addition required to be a morphism of $H$-comodules. In particular, $s^D(1\otimes 1)$ is an $H$-coinvariant.

\begin{thm}
For any comonad $(E,\rho^E,\gamma,D,\rho^D,\Delta^D,\epsilon^D)$ in $\operatorname{ClBimod}^H$ equipped with the $H$-colinear structure morphism $s^D:E\otimes E\to D$, the formula
\begin{equation}
  h\mapsto \gamma(h_{(1)}).s^D(1\otimes 1).\gamma^{-1}(h_{(2)})
\end{equation}
defines a $\genfd$-linear map $H\to D$ which corestricts to a map $y:H\to D^{\mathrm{co}H}$ which is a 1-cocycle in the sense that
\begin{equation}
y(h_{(1)})\otimes_{E^{\mathrm{co}H}} y(h_{(2)}) = \Delta^{\mathrm{co}H}(y(h))
\end{equation}
and $\epsilon^{\mathrm{co}H}(y(h)) = \epsilon(h)1$ 
for all $h\in H$ where $\Delta^{\mathrm{co}H}:D^{\mathrm{co}H}\to D^{\mathrm{co}H}\otimes_{E^{\mathrm{co}H}}D^{\mathrm{co}H}$ is the comultiplication on the $E^{\mathrm{co}H}$-coring of $H$-coinvariants $C=D^{\mathrm{co}H}$ induced by $\Delta^D$, as in Theorem~\ref{thm:coinvcor}. Map $s^D$ admits restriction-corestriction to a map $s^C:U\otimes U\to C$ which is a morphism of $U$-corings.
The equality
\begin{equation}\label{eq:Cmeasfromcocy}
h\triangleright_C s^C(\sum_\alpha u_\alpha\otimes u'_\alpha) = \sum_\alpha(h_{(1)}\triangleright_U u_\alpha).y(h_{(2)}).(h_{(3)}\triangleright_U u'_\alpha)
\end{equation}
holds in $C$ for all $\sum_\alpha u_\alpha\otimes u'_\alpha\in U\otimes U$ and $h\in H$. In particular, if $s^C(u\otimes 1) = s^C(1\otimes v)$ then also
\begin{equation}\label{eq:Umeasintercocy}
(h_{(1)}\triangleright_U u).y(h_{(2)}) = y(h_{(1)}).(h_{(2)}\triangleright_U v)
\end{equation}
holds in $C$.
\end{thm}\begin{proof}
  Since $g=s^D(1\otimes 1)$ is $H$-coinvariant, $\gamma$ a morphism of $H$-comodules and $\gamma^{-1}$ morphism of $H^{\mathrm{co}}$-comodules, $\rho(\gamma(h_{(1)})g\gamma^{-1}(h_{(2)})) = \gamma(h_{(1)(1)})g\gamma^{-1}(h_{(2)(2)})\otimes h_{(1)(2)}h_{(2)(1)}=\gamma(h_{(1)})g\gamma^{-1}(h_{(2)})\otimes 1_H$, hence $y$ corestricts to $D^{\mathrm{co}H}$. The induced $\Delta^{\mathrm{co}H}$ is obtained by restricting $\Delta$ to coinvariants, corestricting to $(D\otimes_E D)^{\mathrm{co}H}$ and postcomposing by isomorphism $(D\otimes_E D)^{\mathrm{co}H}\cong D^{\mathrm{co}H}\otimes_{E^{\mathrm{co}H}} D^{\mathrm{co}H}$. Using this identification, it is sufficient to prove that $\Delta^D(y(h)) = y(h_{(1)})\otimes_E y(h_{(2)})$ holds as an equality in $D\otimes_E D$. Comultiplication $\Delta^D$ is a map of $E$-bimodules so $\Delta^D(\gamma(h_{(1)})g\gamma^{-1}(h_{(2)}))=\gamma(h_{(1)})g\otimes_E g\gamma^{-1}(h_{(2)}) = \gamma(h_{(1)})g\gamma^{-1}(h_{(2)})\gamma(h_{(3)})\otimes_E g\gamma^{-1}(h_{(4)}) = \gamma(h_{(1)})g\gamma^{-1}(h_{(2)})\otimes_E\gamma(h_{(3)})g\gamma^{-1}(h_{(4)}) = y(h_{(1)})\otimes_E y(h_{(2)})$.
\end{proof}

Conversely, given a comonad $(U,\triangleright_U,C,\sigma,\triangleright_C,\Delta^C,\epsilon^C)$ in $\operatorname{Cocyc}^H$, with structure map $s^C:U\otimes U\to C$ and a 1-cocycle $y:H\to C$ we can reconstruct up to an isomorphism a comonad in $\operatorname{ClBimod}^H$ equipped with the $H$-equivariant structure morphism $s^D:E\otimes E\to D$ where $E = U\sharp_{\triangleright_U,\sigma}H$ and $D=C\sharp_{\triangleright_C} H$.

Special cases of covers by localizations \cite{ros:NcSch,skoda:ban,skoda:cech,skoda:gmj,skoda:nloc} and covers by ideals~\cite{brzwri} will be described elsewhere.

\subsubsection{Associated fiber bundles}

Given a left $H$-comodule $F$ (the typical fiber), and a comonad $(E,D,\gamma,s^D)$ in the bicategory $\operatorname{ClBimod}^H$ with $H$-coinvariants coring $(U,C,s^C)$, one can construct an equalizer
$$
\begin{tikzcd}
\Gamma\xi_F\ar[r]& U\otimes F \ar[r,shift left=.75ex,"\kappa_\gamma"]
  \ar[r,shift right=.75ex,swap,"\kappa_\gamma'"]
  &D\square^H F
\end{tikzcd}
$$
The upper arrow is $\kappa_\gamma : u\otimes f\mapsto s^D(u\gamma(f_{(-1)})\otimes 1)\otimes f_{(0)}$ and the lower arrow is $\kappa_\gamma':u\otimes f\mapsto s^D(1\otimes u\gamma(f_{(-1)}))\otimes f_{(0)}$. In terms of the grouplike element $g=s^D(1\otimes 1)$ these maps are respectively $u\otimes f\mapsto u \gamma(f_{(-1)})g\otimes f_{(0)}$ and $u\otimes f\mapsto g u\gamma(f_{(-1)})\otimes f_{(0)}$.

Given the cleaving $\gamma$, the map $e\otimes f\mapsto e\gamma(f_{(-1)})\otimes f_{(0)}$ is a $\genfd$-automorphism of $E\otimes F$ which restricts to isomorphism $U\otimes F\to E\square^H F$  with inverse
$e\otimes f\mapsto e\gamma^{-1}(f_{(-1)})\otimes f_{(0)}$. Using this identification we form an alternative but equivalent equalizer sequence
$$\begin{tikzcd}
\Gamma\xi_F\ar[r]& E\square^H F \ar[r,shift left=.75ex,"\iota_g"]
  \ar[r,shift right=.75ex,swap,"\iota_g'"]
  &D\square^H F
\end{tikzcd}$$
  where $\iota_g:e\otimes f\mapsto e g\otimes f$ and $\iota_g':e\otimes f\mapsto
  g e \otimes f$. This sequence looks more invariant, however the local coordinate patches involve total space, rather than just a base. Now, for all
  $\sum_i u_i\otimes f_i$
  $$
  \sum_i u_i \gamma(f_{i(-1)}) g \otimes f_{i(0)} = \sum_i g u_i f_{i(-1)}\otimes f_{i(0)}
  $$
is an equality in $D\otimes F$ iff
  \begin{equation}
\sum_i u_i y(f_{i(-1)}) \otimes f_{i(0)} = \sum_i g u_i\otimes f_i
  \end{equation}
  
Now, if $F$ is replaced by a left $H$-comodule algebra $\mathcal{F}$ (thus not any more interpreted as a typical fiber but rather function algebra on a typical fiber), then instead we define not a vector space of sections, but, more fundamentally, a noncommutative space $\xi_{\mathcal{F}}$ represented by a usual coring.
Namely, $E\square^H \mathcal{F}$ is then not only a submodule but a subalgebra of the tensor product algebra $E\otimes_\genfd\mathcal{F}$ and
$D\square^H\mathcal{F}$ is an $E\square^H\mathcal{F}$-subbimodule of $D\otimes_\genfd\mathcal{F}$. Given any Hopf modules $P\in {}_E\mathcal{M}^H$ the canonical map
$$
c_{P,\mathcal{F}}:D\square^H\mathcal{F}\otimes_{E\square^H\mathcal{F}}P\square^H\mathcal{F}
\to (D\otimes_E P)\square^H\mathcal{F}
$$
given by $(d\otimes f)\otimes(p\otimes f')\mapsto d\otimes p\otimes f f'$ is an isomorphism provided $E$ is faithfully flat $H$-Hopf-Galois.
Specialize to $P=D$ to define the new coproduct 
$$
\Delta_{\xi_{\mathcal{F}}} = c^{-1}_{D,\mathcal{F}}\circ(\Delta^D\square^H\id_{\mathcal{F}}):D\square^H\mathcal{F}\to D\square^H\mathcal{F}\otimes_{E\square^H\mathcal{F}} D\square^H\mathcal{F}
$$
One proves~\cite{skoda:ass} that $(D\square^H\mathcal{F}, \Delta_{\xi_{\mathcal{F}}},\epsilon_{\xi_{\mathcal{F}}})$ is an $E\square^H\mathcal{F}$-coring, which we view as a representation of a total space of an associated bundle $\xi_{\mathcal{F}}$; it comes along with a morphism from the base $E^{\mathrm{co}H}$-coring $D^{\mathrm{co}H}$ dualizing the projection of the associated bundle $\xi_{\mathcal{F}}$.

\end{document}